\documentclass{article}

\usepackage{amsmath}
\usepackage{amsthm}
\usepackage{amssymb}
\usepackage{amsfonts}
\usepackage{latexsym}
\usepackage{delarray}
\usepackage[dvips]{graphics}
\usepackage{epsfig}
\usepackage{color}

\newtheorem{thm}{Theorem}[section]
\newtheorem*{thm*}{Theorem}
\newtheorem{lemma}[thm]{Lemma}

\newtheorem{cor}[thm]{Corollary}
\newtheorem{question}{Question}

\theoremstyle{definition}

\theoremstyle{remark}
\newtheorem*{obs}{Observation}

\newcommand{\Aut}{\operatorname{Aut}}

\title{The maximum distinguishing number of a group}
\author{Melody Chan \\ University of Cambridge \\ Cambridge, England \\ \texttt{melody.chan@aya.yale.edu}}

\begin{document}

   \date{
   \small Mathematics Subject Classification: 05E15, 20B25, 20D60}
   
\maketitle

\begin{abstract}
Let $G$ be a group acting faithfully on a set $X$.  The distinguishing number of the action of $G$ on $X$, denoted $D_G(X)$, is the smallest number of colors such that there exists a coloring of $X$ where no nontrivial group element induces a color-preserving permutation of $X$.  In this paper, 
we show that if $G$ is nilpotent of class $c$ or supersolvable of length $c$ then $G$ always acts with distinguishing number at most $c+1$.  We obtain that all metacyclic groups act with distinguishing number at most 3; these include all groups of squarefree order.  We also prove that the distinguishing number of the action of the general linear group $GL_n(K)$ over a field $K$ on the vector space $K^n$ is 2 if $K$ has at least $n+1$ elements.
\end{abstract}

\section{Introduction}

An action of a group $G$ on a set $X$ is said to be {\em faithful} if only the identity element of $G$ fixes every element of $X$.  Let $G$ be a group acting faithfully on $X$.  For $r \in \mathbb{N}$, an {\em $r$-coloring} of $X$ is  a function $c \colon  X \rightarrow \{1, \ldots, r\}$.  A permutation $\pi$ of $X$ {\em preserves} the coloring $c$ if $c(\pi x) = c(x)$ for all $x \in X$.  A coloring is said to be {\em distinguishing} if the only element in $G$ that induces a color-preserving permutation of $X$ is the identity element.  The {\em distinguishing number} of the action of $G$ on $X$, denoted $D_G(X)$, is the smallest $r$ admitting a distinguishing $r$-coloring of $X$ with respect to the action of $G$.  If there does not exist a distinguishing $r$-coloring of $X$ for any finite $r$, we say that $D_G(X) = \infty$.  

The distinguishing number was first defined by Albertson and Collins in \cite{a_c} as a property of graphs.  More specifically, the distinguishing number of a graph $M$, denoted $D(M)$, is the smallest number of colors admitting a coloring of the vertices such that the only color-preserving automorphism of $M$ is the identity; thus $D(M) = D_{\Aut(M)}(V(M))$.  Note that distinguishing colorings of graphs need not be proper colorings in the graph theoretic sense: two adjacent vertices may or may not have the same color.  Although Albertson and Collins initially defined the distinguishing number solely in terms of graphs, the approach they chose to take is nevertheless highly group theoretic.  Given a group $G$, they define the \emph{distinguishing set} of $G$, denoted $D(G)$, as

$$D(G) = \{ D(\mathcal{G})~|~\mathcal{G}\textrm{ is a graph with }\Aut(\mathcal{G}) \cong G \}.$$

Their results in \cite{a_c} center around characterizing the distinguishing set of a group.  For example, they show that $D(S_4) = \{2,4\}$.  They also prove the following result.

\begin{thm*} \cite[Corollary 3.1, Theorem 6]{a_c}

(1) If $G$ is Abelian then $\max\{D(G)\} \le 2$.

(2) If $G$ is dihedral then $\max\{D(G)\} \le 3$.
\end{thm*}

In addition, the distinguishing number of several families of graphs, including trees, hypercubes, and generalized Petersen graphs, have been studied in \cite{b_c}, \cite{mc_cubes}, \cite{cctcheng}, and \cite{potanka}.

In \cite{jt}, Tymoczko generalizes the notion of the distinguishing number to group actions on sets and proves results about the distinguishing number of actions of the symmetric group $S_n$.
She shows that the distinguishing number of group actions is indeed a more general question than the distinguishing number of graphs.  For example, she exhibits a faithful $S_4$-action with distinguishing number 3, contrasting Albertson and Collins' result that $D(S_4) = \{2,4\}$.  This difference highlights the fact that not all faithful group actions are realized as actions of the automorphism group of a graph on its vertex set.

Following Tymoczko, it seems natural to expand the notion of the distinguishing set of a group to include all of its possible actions, not just those arising from graph automorphism groups.  In this generalized context, we ask the following question: given a group $G$, what is the best upper bound we can give for $D_G(X)$?  In Section \ref{maximum}, we give an upper bound for the maximum distinguishing number for a large class of groups including nilpotent and supersolvable groups.  

\begin{thm*} If $G$ is nilpotent of class $c$ or supersolvable of length $c$ then $G$ acts with distinguishing number at most $c+1$.  
\end{thm*}

As a corollary, we obtain that all metacyclic groups act with distinguishing number at most 3 since they are supersolvable of length 2; these include all groups of squarefree order.  Albertson and Collins' results for Abelian and dihedral groups are special cases of nilpotent groups of class 1 and metacyclic groups, respectively.  In Section \ref{generallinear}, we compute the distinguishing number for an important group action, that of the general linear group over a field $K$ on a vector space over $K$.  We show that the distinguishing number of this action is 2 if $|K| > n+1$ where $n$ denotes the dimension of the vector space.

Our definition of the distinguishing number of a group action differs from the one given in \cite{jt} in that we require the action to be faithful.  This apparent restriction does not in actuality limit the question being considered, however, for given a nonfaithful action of $G$ on $X$, we may consider instead the faithful action of the quotient group $G/Stab(X)$ on $X$, where $Stab(X)$ denotes the elements of $G$ that fix each $x \in X$.  Also, in contrast to both \cite{a_c} and \cite{jt}, we do not require our groups and sets to be finite, simply because there seems to be no reason to do so.  We only note that if $G$ is an infinite group acting faithfully on a set $X$, then $X$ must be infinite as well.

\section{The maximum distinguishing number of a group} \label{maximum}

Given a group $G$, let $\overline{D}(G) = \max \bigl \{D_G(X)~|~G \textrm{ acts faithfully on } X \bigr \}$ denote the largest distinguishing number that $G$ admits, or $\overline{D}(G) = \infty$ if some $D_G(X) = \infty$.  In \cite{jt}, Tymoczko proves the bound $D_G(X) \le k$ if $|G| \le k!$, a result originally formulated by Albertson, Collins and Kleitman in terms of graphs.  This result holds for any action of $G$ on $X$, faithful or nonfaithful.  In \cite{a_c}, Albertson and Collins show that if $G$ is Abelian then $\overline{D}(G) \le 2$ and if $G$ is dihedral then $\overline{D}(G) \le 3$.  Their proof is formulated only in terms of graphs but also holds for group actions.  In this section, we generalize these results to a class of groups that includes all nilpotent and supersolvable groups.

The following lemma gives some conditions under which we may characterize the maximum distinguishing number of a group.  The idea to consider the intersection of a normal subgroup with the stabilizing subgroup of orbit representatives was inspired by Albertson and Collins' proof for dihedral groups in \cite{a_c}.  In what follows, we use $\langle x \rangle$ to denote the subgroup generated by a group element $x$.  Also, we will denote group actions by exponentiation on the right.  Thus, the image of an element $x \in X$ under the action of $g \in G$ is denoted $x^g$, and we have $(x^{g_1})^{g_2} = x^{(g_1 g_2)}$ for all $g_1, g_2 \in G$.  

\begin{lemma} \label{l:extension}
Suppose $N$ is a normal subgroup of $G$ with the property that if $n_1, n_2 \in N$ are conjugate elements in $G$, then $\langle n_1 \rangle = \langle n_2 \rangle $.  Suppose further that any subgroup $L$ of $G/N$ has the property that $\overline{D}(L) \le c$.  Then $\overline{D}(G) \le c+1$.
\end{lemma}

\begin{proof}  The case $G = 1$ is trivial.  Suppose that a nontrivial group $G$ acts faithfully on a set $X$.  Choose a set $U$ of representatives of the orbits of $G$ on $X$ (using the Axiom of Choice if there are infinitely many orbits), and let $H = \{g \in G~|~u^g = u \textrm{ for each } u \in U\}$ stabilize the set $U$ pointwise.  We claim that $H \cap N = 1$.
Suppose that $n \in H \cap N$, so that $n$ stabilizes each $u \in U$.  Fix any $x \in X$ and let $u \in U$ be the representative of the orbit containing $x$.  Let $g \in G$ satisfy $u = x^g$, and let $H_x$ be the stabilizer subgroup of $x$.  By assumption, $n$ stabilizes $u$, so $g n g^{-1}$ stabilizes $x$.  But the fact that $\langle n \rangle = \langle g n g^{-1} \rangle$ implies that $n \in \langle g n g^{-1} \rangle \le H_x$.  Therefore, $n$ stabilizes each $x \in X$.  Since $G$ acts faithfully, $n = 1$ and so $H \cap N = 1$.  Applying the Second Isomorphism Theorem, we conclude that $H \cong HN/N \le G/N$.  

Now, we know that $X \setminus U$ is nonempty because $G$ is nontrivial, so consider the action of $H$ on $X \setminus U$.  This action is faithful since the action of $G$ on $X$ is faithful.  Then $D_H(X \setminus U) \le \overline{D}(H) \le c$ since $H$ is isomorphic to a subgroup of $G/N$.  Then let $C\colon  X \setminus U \rightarrow \{1, \ldots, c \}$ be a $c$-coloring of $X \setminus U$ that is distinguishing with respect to the action of $H$.  Now define $C'\colon  X \rightarrow \{1, \ldots, c+1\}$ as 
$$
C'(x) =
\begin{cases}
c+1 & \textrm { if } x \in U, \\
C(x) & \textrm{ if } x \not \in U.
\end{cases}
$$
We claim that $C'$ is a distinguishing $(c+1)$-coloring of $X$ with respect to the action of $G$.  Suppose $g \in G$ preserves $C'$.  Then $g$ must fix each orbit representative $u \in U$, since they are the only elements of color $c+1$ and lie in different orbits.  Thus $g \in H$.  Then consider the action of $g$ on $X \setminus U$.  Since the restriction of $C'$ to $X \setminus U$ is a distinguishing coloring with respect to the action of $H$, and $g \in H$ preserves this coloring, we have $g = 1$.  Therefore, $C'$ is a distinguishing $c+1$-coloring with respect to the action of $G$.  We conclude that an arbitrary action of $G$ has distinguishing number at most $c+1$, and so $\overline{D}(G) \le c+1$.

\end{proof}

The next theorem is a consequence of Lemma \ref{l:extension}.  Following \cite{hall}, we define a \emph{normal series} for a group $G$ to be a chain of subgroups $1 = G_0 \triangleleft G_1 \triangleleft \cdots \triangleleft G_c = G$ with the additional condition that each $G_i \triangleleft G$.  

\begin{thm} \label{t:normalseries}
Suppose a group $G$ has a finite normal series
$$1 = G_0 \triangleleft G_1 \triangleleft \cdots \triangleleft G_c = G$$
in which each quotient $G_{i+1}/G_i$ is cyclic or is contained in $Z(G/G_i)$.  Then $\overline{D}(G) \le c+1$.
\end{thm}
\begin{proof}
We proceed by induction on $c$.  If $c = 0$, then $G = G_0 = 1$ and $\overline{D}(G) = 1$.  Now let $G$ have a normal series $1 = G_0 \triangleleft G_1 \triangleleft \cdots \triangleleft G_c = G$ of length $c > 0$ with the required property.  In order to apply Lemma \ref{l:extension}, we wish to show that any two conjugate elements lying in $G_1$ generate the same subgroup, and in addition any subgroup of the quotient group $G/G_1$ acts with distinguishing number at most $c$.

Let $n_1$ and $n_2$ be conjugate elements in $G_1$.  We have assumed that $G_1$ is either cyclic or contained in $Z(G)$.  In the former case, note that since conjugation is a group automorphism, $[G_1:\langle n_1 \rangle] = [G_1:\langle n_2 \rangle]$.  But $G_1$ is cyclic, so it has precisely one subgroup of this index.   Therefore $\langle n_1 \rangle = \langle n_2 \rangle$.  In the latter case, note that every element of $Z(G)$ has no conjugates other than itself, so $n_1 = n_2$ and $\langle n_1 \rangle = \langle n_2 \rangle$.

Next, it follows from the Third Isomorphism Theorem that $$1 = G_1/G_1 \triangleleft G_2/G_1 \triangleleft \cdots \triangleleft G/G_1$$ is a normal series for $G/G_1$ of length $c-1$ in which each quotient group $(G_{i+1}/G_1)/(G_i/G_1))$ is cyclic or is contained in $Z((G/G_1)/(G_i/G_1))$.  Now for any subgroup $L$ of $G/G_1$, let $L_i = (G_{i+1}/G_1) \cap L$ for each $i$ with $0 \le i \le c-1$.  Then one can check that 
$$1 = L_0 \triangleleft L_1 \triangleleft \cdots \triangleleft L_{c-1} = L$$
is a normal series of length $c-1$ for $L$ with the property that each quotient group $L_{i+1}/L_i$ is cyclic or is contained in $Z(L/L_i)$.  Then $\overline{D}(L) \le c$ by the inductive hypothesis.  Thus, all the conditions of Lemma \ref{l:extension} are satisfied, so $\overline{D}(G) \le c+1$.
\end{proof}

As consequences of Theorem \ref{t:normalseries}, we obtain upper bounds on the distinguishing number of nilpotent and supersolvable groups.  We recall the definitions of these important classes of groups below; see \cite{hall} for a more detailed discussion of them. 

A group $G$ is said to be \emph{nilpotent} if it possesses a finite normal series $1 = G_0 \triangleleft G_1 \triangleleft \cdots \triangleleft G_c = G$ such that each quotient group $G_{i+1}/G_i$ is contained in $Z(G/G_i)$.  If the shortest such normal series has length $c$, then we say that $G$ is \emph{nilpotent of class $c$}.

\begin{cor} \label{c:nilpotent}
Let $G$ be nilpotent of class $c$.  Then $\overline{D}(G) \le c+1$.
\end{cor}

In particular, since all Abelian groups are class-1 nilpotent, we have $\overline{D}(G) \le 2$ for $G$ Abelian, as shown in \cite{a_c}.

A group $G$ is said to be \emph{supersolvable} if it possesses a finite normal series $1 = G_0 \triangleleft G_1 \triangleleft \cdots \triangleleft G_c = G$ such that each quotient group $G_{i+1}/G_i$ is cyclic.  In this case, we will say that $G$ is \emph{supersolvable of length $c$}.  See \cite{bray} for a detailed discussion of supersolvable groups.

\begin{cor} \label{c:supersolvable}
Let $G$ be supersolvable of length $c$.  Then $\overline{D}(G) \le c+1$.
\end{cor} 

A group $G$ is called {\em metacyclic} if it has a normal subgroup $N \triangleleft G$ such that both $N$ and $G/N$ are cyclic.  Such groups have been completely classified in \cite{hempel}, and include all groups of squarefree order.

\begin{cor}
Let $G$ be a metacyclic group.  Then $\overline{D}(G) \le 3$.
\end{cor}

We obtain as a special case that if $G$ is dihedral then $\overline{D}(G) \le 3$, as shown in \cite{a_c}.

\section{The action of $GL_n(K)$ on $K^n$} \label{generallinear}

In this section, we consider the action of $GL_n(K)$, the group of $n \times n$ invertible matrices over a field $K$, on $K^n$, the $n$-dimensional vector space over $K$.  We may regard the elements of $K^n$ as column vectors and accordingly define a left action of $GL_n(K)$ on $K^n$ as $A\colon  v \mapsto Av$ for each $v \in K^n$, $A \in GL_n(K)$.  This action is clearly faithful.  

Our main result is that if $K$ is sufficiently large, then 2 colors suffice to distinguish this action.

\begin{thm} \label{t:glnk}
Let $K$ be a field.  If $K$ is infinite or is finite of order greater than $n+1$, then $D_{GL_n(K)}(K^n) = 2$.
\end{thm}

\begin{proof}
We first observe that the multiplicative group $K^\times$ must contain a nonzero element $\alpha$ of order greater than $n$.  For if $K$ is infinite, then we may certainly choose such an $\alpha$ since there exist only finitely many solutions in $K$ to the equations $x^l = 1$ for each $1 \le l \le n$.  On the other hand, if $K$ is finite, then we know that $K^\times$ is a cyclic group of order $|K| - 1$.  Let $\alpha$ generate the group $K^\times$, then the order of $\alpha$ is $|K| - 1 > n$.

Now let $e_1, \ldots, e_n$ be the standard basis vectors in $K^n$, and let $S$ be the set of vectors $\{ \alpha^i e_j~|~0 \le i < j \le n \}$.  Each of these vectors is distinct since $\alpha$ has order greater than $n$, so the cardinality of $S$ is $\frac{1}{2}n(n+1)$.  Now color every vector in $S$ blue and all remaining vectors red.  We claim this is a distinguishing 2-coloring of $K^n$ with respect to the action of $GL_n(K)$.  

Suppose $A \in GL_n(K)$ preserves this coloring.  It suffices to show that $Ae_k = e_k$ for each $e_k$.  Since $e_k$ is blue, the image of $e_k$ must also be blue and so must have the form $\alpha^i e_j$ for $0 \le i < j \le n$.  We wish to show that $i = 0$ and $j = k$.  First, note that $Ae_k = \alpha^i e_j$ implies that $A^{-1}(\alpha^{i-1}e_j) = \alpha^{-1} A^{-1}(\alpha^i e_j) = \alpha^{-1} e_k$.  Now, $\alpha^{-1}e_k$ is a red point, because if instead $\alpha^{-1}e_k = \alpha^c e_k$ for some $0 \le c < k$, then $\alpha$ would have order at most $c+1$, but $c+1 \le k \le n$ and we assumed that the order of $\alpha$ was greater than $n$.  So, $\alpha^{-1}e_k$ is red, and since $A$ was assumed to be color preserving, $\alpha^{i-1}e_j$ is also red.  This is only possible if $i = 0$.  Thus $Ae_k = e_j$ and so $A$ induces some permutation of the basis vectors $\{ e_1, \ldots, e_n \}$.  Suppose for a contradiction that $A$ permutes them nontrivially.  Then there must exist $Ae_k = e_j$ with $k < j$.  Then $\alpha^{j-1}e_k$ is red since $k \le j-1$, but $\alpha^{j-1}e_j$ is blue, and $A(\alpha^{j-1}e_k) = \alpha^{j-1}e_j$, a contradiction.  Therefore, $A$ fixes each basis vector $e_k$ and so $A = \mathbf{1}_n$ as desired.  We have exhibited a distinguishing $2$-coloring of $K^n$, so $D_{GL_n(K)}(K^n) \leq 2$.  Now, it is possible that $D_{GL_n(K)}(K^n) = 1$ only if $GL_n(K)$ is the trivial group.  This occurs only when $n = 1$ and $K = \mathbb{F}_2$, which was excluded by the assumption that $|K| > n+1$.  Therefore, we have the equality $D_{GL_n(K)}(K^n) = 2$.

\end{proof}

Theorem \ref{t:glnk} leaves open the case when the size of the field is relatively small.  It is possible to show by case analysis that $D_{GL_2(\mathbb{F}_2)}(\mathbb{F}_2^2) = D_{GL_2(\mathbb{F}_3)}(\mathbb{F}_3^2) = 3$.  However, we leave the case when $n > 2$ and $|K| \le n+1$ as an open problem.

\section{Discussion and open questions}
The distinguishing number seems to be a very natural property of group actions, and efforts to relate the distinguishing number of a group action to group properties seem likely to be fruitful.  Below, we give several possibilities for further investigation.

One interesting method of attack relies on the following simple fact.

\begin{obs} Let $G$ act faithfully on $X$. Fix a coloring $c$ of $X$ and let $H_c = \{ g \in G ~|~ g \textrm{ preserves } c \}$.  Then $H_c$ is a subgroup of $G$.
\end{obs}

The distinguishing number, then, is the smallest number of colors admitting a coloring $c$ such that $H_c = 1$.  Thus, it seems plausible that one could make direct use of information on the subgroup structure of $G$ to characterize the distinguishing number.  We present the theorem below as an example of employing this technique.

\begin{thm}
Let $G$ be a finite group acting faithfully on a set $X$.  Let $p$ be the smallest prime dividing the order of $G$, and let $M$ be the length of the largest orbit of the action of $G$ on $X$.  Then $D_G(X) \le \lceil \frac{M}{p-1} \rceil$.
\end{thm}

\begin{proof}  With $\lceil \frac{M}{p-1} \rceil$ colors, we may color each orbit of $X$ such that every color class within a given orbit has size at most $p-1$.  Call this coloring $c$.  Let $\sim$ be the equivalence relation given by $x_1 \sim x_2$ if and only if $x_1$ and $x_2$ are in the same orbit and have the same color. Let $P_1, P_2, \ldots, P_k$ be the equivalence classes of this relation, and let $n_i = |P_i|$.  Note that a color preserving permutation $h \in H_c$ can take a given element only to another element in its equivalence class.  Thus $H_c \le S_{n_1} \times S_{n_2} \times \cdots \times S_{n_k}$.  Also, $H_c \le G$.  But since each $n_i < p$, the orders of $S_{n_1} \times S_{n_2} \times \cdots \times S_{n_k}$ and $G$ are relatively prime.  Since the order of $H_c$ divides both orders, we have $H_c = 1$, and so $c$ is distinguishing.
\end{proof}

In addition, we could use subgroup structure as a way to generalize the notion of the distinguishing number, as follows.  Given a group $G$ acting faithfully on a set $X$ and $H$ a subgroup of $G$, let $D_{G,H}(X)$ denote the smallest number of colors admitting a coloring of $X$ such that the only elements of $G$ that induce color-preserving permutations lie in $H$.  Thus, when $H = 1$, we recover the original notion of the distinguishing number.

\begin{question}
Characterize $D_{G,H}(X)$.
\end{question}

In Section \ref{maximum}, we considered the maximum distinguishing number admitted by a given group.  Intuitively, we would expect a large group to admit actions that require many colors to distinguish them.  Thus, we ask whether the distinguishing number is ordered in a way that respects the partial ordering of groups defined by subgroup inclusion.  
\begin{question}
Let $G$ and $H$ be groups, $H$ a subgroup of $G$.
Does it follow that $\overline{D}(H) \le \overline{D}(G)$?
\end{question}
Note that if a given faithful action of $H$ on $X$ can be extended to a faithful action of $G$ on $X$ then $D_H(X) \le D_G(X)$, because any coloring of $X$ that is distinguishing with respect to the action of $G$ is also distinguishing with respect to the action of $H$.  However, since not every faithful action of $H$ on $X$ can necessarily be extended to a faithful action of $G$ (for example if $|G| > |X|!$), the question cannot be answered immediately in the affirmative.

We also ask whether the bounds obtained in Section \ref{maximum} for nilpotent and supersolvable groups are tight.

\begin{question}
For which $k > 2$ does there exist a group $G$ that is nilpotent group of class $k$ (or a supersolvable group of length $k$) acting faithfully on a set $X$ such that $D_G(X) = k+1$?
\end{question}

In Section \ref{generallinear}, we showed that $D_{GL_n(K)}(K^n) = 2$ if $|K| > n+1$ and $D_{GL_n(K)}(K^n) = 3$ if $|K| \in \{2,3\}$ and $n = 2$.  As mentioned, we leave the remaining cases as an open question.

\begin{question}
Compute $D_{GL_n(K)}(K^n)$ for $n > 2$ and $|K| \le n+1$.
\end{question}

The generalization of the distinguishing number to infinite groups acting on infinite sets is new, and it might be interesting to investigate conditions on the finiteness or infiniteness of the distinguishing number.  This leads to many questions, including the following.

\begin{question}
Suppose $G$ is a group that always acts with finite distinguishing number.  Does it follow that $\overline{D}(G) < \infty$, that is, that the set $\{D_G(X)\}$ has a maximum element?
\end{question}

Another approach would simply be to define the distinguishing number of an action to be the cardinality of the smallest set of colors admitting a distinguishing coloring.  This would eliminate the formal distinction between the finite and infinite cases.

Finally, throughout this paper, we have considered the distinguishing numbers that a fixed group admits in its actions on various sets.  We could instead fix a set $[n]$ and consider the distinguishing numbers it admits under the actions of various groups.
\begin{question}
For each $n$, characterize the set $$T_n = \{D_G([n])~|~G \textrm{ is a transitive subgroup of } S_n\}.$$
\end{question}
One may show that $T_n = \{2, \ldots, n\}$ for $n = 2,3,4,5,$ and $6$.  We ask whether $T_n$ has this form for larger $n$.  Note that we require our group $G$ to be transitive, for otherwise each distinguishing number $k$ between $1$ and $n$ could be achieved by taking a subgroup of $S_n$ that fixes each $k+1, k+2, \ldots, n$ and whose action on $1, \ldots, k$ is isomorphic to the action of $S_k$.

\section{Acknowledgments}
This research was conducted at the University of Minnesota Duluth Research Experience for Undergraduates, while the author was a student at Yale University.  The author would like to express her thanks to Melanie Wood for numerous ideas and suggestions on drafts of this paper, to Philip Matchett and Daniel Isaksen for several very helpful conversations, and to Joseph Gallian for his support.  This research was funded by the National Science Foundation (DMS-0137611) and the National Security Agency (H-98230-04-1-0050).

\end{document}